\documentclass{article}
\usepackage{a4}
\usepackage{amssymb}
\usepackage{epsfig}
\newcommand{\sect}[1]{\setcounter{equation}{0}\section{#1}}

\renewcommand{\appendix}{\setcounter{section}{0}
\renewcommand{\thesection}{\Alph{section}}}

\def\hR{\hat R}
\def\hx{\hat x}
\def\bdm{\begin{displaymath}}
\def\edm{\end{displaymath}}
\def\be{\begin{equation}}
\def\ee{\end{equation}}
\def\ba{\begin{eqnarray}}
\def\ea{\end{eqnarray}}
\def\T{T_{orb}}
\newcommand{\nn}{\nonumber}
\newcommand{\no}{\nonumber \\}

\def\ha{\frac{1}{2}}
\def\l{\lambda}
\def\b#1{{\mathbb #1}}
\def\hq{\frac{1}{q}}
\def\nq{\frac{1}{q^2}}
\def\pn{\par\noindent}

\def\ket#1{\:\:\vert #1 \rangle}
\begin{document}
\begin{titlepage}
\rightline{LMU-TPW 2000-06}
\rightline{MPI-PhT/2000-08}
\vspace{4em}
\begin{center}
 
{\Large{\bf Structure of the Three-dimensional \\ Quantum Euclidean Space}}
 
\vskip 3em
 
{{\bf
B.L. \ Cerchiai${}^{1,2}$, J.\ Madore${}^{2,3}$, S.\ Schraml${}^{1,2}$,
J.\ Wess${}^{1,2}$ }}
 
\vskip 2em

${}^{1}$Sektion Physik,
Universit\"at M\"unchen\\
Theresienstr.\ 37,
D-80333 M\"unchen \\[1em]
${}^{2}$Max-Planck-Institut f\"ur Physik\\
        F\"ohringer Ring 6, D-80805 M\"unchen\\[1em]
${}^3$Laboratoire de Physique Th\'eorique\\
        Universit\'e de Paris-Sud, B\^atiment 211, F-91405 Orsay\\[1em]
\end{center}
 
\vspace{2em}   
\begin{abstract}
As an example of a noncommutative space we discuss the quantum
3-dimensional Euclidean space $\b R^3_q$ together with its 
symmetry structure in great detail.
The algebraic structure and the representation theory 
are clarified and discrete spectra for the coordinates are found.
The $q$-deformed Legendre functions play a special role. A completeness
relation is derived for these functions.
\end{abstract}
\end{titlepage}

\newpage
\sect{Introduction}

Noncommutative space-time structures are finding increasing interest in gauge
theories. Special examples \cite{alle} have been discussed in the literature. Among
them are 
\medskip
\pn
1) The canonical structure:
\bdm 
[\hat{x}^i, \hat{x}^j] = i \Theta^{ij}, 
\qquad \Theta^{ij} \, \in \, \mathbb{C}
\edm
2) the Lie-algebra structure:
\bdm
[\hat{x}^i, \hat{x}^j] = i\, c^{ij}_k \hat{x}^k, 
        \qquad c^{ij}_k \, \in \, \mathbb{C}
\edm
3) The quantum space structure:
\bdm
[\hx^i, \hx^j] = i c^{ij}_{kl} \hx^k \hx^l, 
\qquad c^{ij}_{kl} \, \in \, \mathbb{C}
\edm

We shall discuss a special example of the third case~\cite{FadResTak89}. 
For a more general review of the quantum space structure see \cite{Wess}. 

In all of the above cases we consider the associative algebra freely 
generated by the elements $\hx^i$ modulo the respective relations.
This algebra of formal power series forms the algebra 
$\mathcal{A}_x$.
\bdm
\mathcal{A}_x = \frac{\mathbb{C}[\,[\hx^1 \ldots \hx^N]\,]}{R}
\edm
For a physicist this means that he is free to use the relations to 
reorder the elements of an arbitary power series. For the quantum space
algebra we wish to exclude pathological cases such as the trivial case
where there is no relation at all or where the product of any two
elements is zero modulo the relations. To exclude such cases we shall demand
the Poincar\'{e}-Birkhoff-Witt property for the algebra. By this we mean
that the dimension of the space of homogeneous polynomials is the
same as in the commutative case. For the first and second examples this 
will be the case, for the third example we require the Yang-Baxter
equation. To formulate it we write the relations in the form:
\bdm
\hx^i \hx^j = \hR^{ij}_{kl} \hx^k \hx^l, \qquad \hR^{ij}_{kl} \, \in \,
        \mathbb{C}
\edm
and define $N^3 \times N^3$ matrices
\ba
\hR^{\quad i_1 \, i_2 \, i_3}_{12 \, j_1 \, j_2 \, j_3} &=& 
\hR^{i_1 \, i_2}_{j_1 \, j_2}\, \delta^{i_3}_{j_3} \nn \\
\hR^{\quad i_1 \, i_2 \, i_3}_{23 \, j_1 \, j_2 \, j_3} &=& 
\delta^{i_1}_{j_1}\, \hR^{i_2 \, i_3}_{j_2 \, j_3} \nn
\ea
The Yang-Baxter equation is:
\bdm
\hR_{12}\,\hR_{23}\,\hR_{12}\,=\,\hR_{23}\,\hR_{12}\,\hR_{23}
\edm
There are several known solutions of this equation. We are interested 
in such relations that allow a conjugation which makes the algebra a
$*$-algebra. This is because we have to associate the observables 
like the coordinates with selfadjoint operators in a Hilbertspace.
The $\hR$-matrices for the quantum groups $SO_q(n)$ allow such conjugations.

The quantum space algebra is a comodule of a quantum group. We start from 
co-algebra relations
\bdm
\Delta (x^i) = T^i_{\, k} \otimes \hx^k
\edm
and compute
\ba
\Delta(x^i) \Delta(x^j) &=& T^i_{\, k}T^j_{\, l} \otimes \hx^k \hx^l \nn \\
&=& T^i_{\, k} T^j_{\, l} \otimes \hR^{kl}_{\, mn} \hx^m \hx^n \nn
\ea
If we demand $RTT$ relations 
\bdm
T^i_{\, k}T^j_{\, l}\hR^{kl}_{\, mn}=\hR^{ij}_{\, kl}T^k_{\, m} T^l_{\, n}
\edm
for the $T$-algebra, we find
\bdm
\Delta(x^i) \Delta(x^j) = \hR^{ij}_{\, kl} \Delta(x^k) \Delta(x^j)
\edm
There is always a solution to the $RTT$ relations given by
\bdm
T^i_{\, j} = \delta^i_j.
\edm
If this is the only solution then the bialgebra consists of the
unit element only; not very interesting. If the $\hR$-matrix leads
to the quantum group $SO_q(n)$, we have a more interesting case. Instead of
introducing the quantum group $SO_q(n)$ we shall deal with the corresponding
$q$-Lie algebra $so_q(n)$. The quantum space is then a module of this
algebra. 

In this paper we discuss the 3-dimensional case in great detail. The
algebra  is introduced in Chapter 2. It has a peculiar property, there is
a homomorphism of the algebra $so_q(3)$ into
the algebra $\b R^3_q$. This is discussed in Chapter 3.
The full algebra can then be generated by a central element, the radius $R$, 
and elements of the tensor product of an $su_q(2)$ algebra and an 
$su_q(1,1)$ algebra. The 
generators of the $su_q(2)$ algebra are further restricted by relations
that when the algebra is represented lead to a unique infinite-dimensional
representation of $su_q(2)$. We call this algebra the $t$-algebra.
The $su_q(1,1)$ algebra we call $K$-algebra. If we then demand that the
$so_q(3)$ algebra corresponds to orbital angular momentum the $K$ algebra
is restricted in the same sense as the $t$-algebra. This is discussed in
Chapter~4. This clarifies the algebraic structure of the $so_q(3)$ module
$\b R^3_q$.

To discuss physics we need representations of the algebra. The observables
should be represented by (essentially) self-adjoint linear operators in a
Hilbert space. This way we can use the well developped formalism of quantum 
mechanics and its interpretation scheme.
In Chapter 5 we discuss the representations of the algebra. We find that
they are characterized by one real parameter $z_0$. In all these
representations we obtain a discrete spectrum for the coordinate $X^3$,
which along with $R$ and $\T^3$, the third component of the orbital angular
momentum, form a complete commuting set of observables. The scale of the
spacing of the eigenvalues of $X^3$ is determined by the constant $z_0$,
the eigenvalues are exponentially spaced. This we call a $q$-lattice.
We are not  surprised that noncommuting variables lead to a discretization
(latticization) of space~\cite{Sny47a}.

In Chapter 6 we construct the transformation that leads to a basis where
$\vec \T^2$ is diagonal. The corresponding transformation function
turn out to be the $q$-deformed associated Legendre functions. They
are defined in Appendix \ref{appD} in terms of the big $q$-Jacobi polynomials.
They satisfy a difference equation, a recursion equation and have
orthogonaltity properties - in complete analogy to the usual associated
Legendre functions. From the self-adjointness property of $X^3$ we
derive a completeness relation as well, this is done in Appendix \ref{appE}. Appendices
\ref{appA}, \ref{appB} and \ref{appC} are devoted to the representation of the $su_q(2)$ and $su_q(1,1)$
algebras and their comultiplication.

\sect{The Algebra of the Euclidean Quantum Space~\boldmath{$\b R^3_q$}}

This algebra has been discussed in Ref \cite{L.W.W.}, we use the same notation
here:

\pn
$\b R^3_q$:
\ba
X^3 X^+ -q^2 X^+ X^3&=&0 \no
X^3 X^- -q^{-2} X^- X^3&=&0                     \label{1.1}\\
X^- X^+ - X^+ X^-&=&\l X^3 X^3,\qquad \lambda=q-q^{-1}, \quad q \in \b{R}. \nn
\ea
We shall assume $q>1$ in this paper. 
This non-commutative structure is our model for a non-commutative
space. We can impose conjugation properties that are compatible with the
relations (\ref{1.1}) justifying  the `$\b{R}$' in $\b{R}^3_q$:
\be
\overline{X^+}=-qX^-, \qquad \overline{X^3}=X^3.               \label{1.2}
\ee

The quantum space $\b{R}^3_q$ has a co-module structure under the action of the 
quantum group $SO_q(3)$ \cite{Wess} and a module structure under
the corresponding $q$-Lie algebra.

\pn
$su_q(2)$:
\ba
q^{-1} T^+ T^--q T^-T^+&=&T^3\no
q^2 T^3 T^+-q^{-2} T^+ T^3&=&(q+q^{-1}) T^+ \label{1.3}\\
q^2 T^- T^3-q^{-2}T^3 T^-&=&(q+q^{-1}) T^- \nn
\ea
The conjugation properties justifying the `$u$' in $su_q(2)$ are:
\be
\overline{T^+}=\frac{1}{q^2}T^-, \qquad \overline{T^3}=T^3. \label{1.4}
\ee

The module structure that was found in Ref. \cite{L.W.W.} is:
\ba
T^3 X^3&=&X^3 T^3 \no
T^3 X^+&=&q^{-4} X^+T^3+q^{-1}(1+q^{-2})X^+ \label{1.5} \\
T^3 X^-&=&q^4 X^-T^3-q(1+q^2)X^- \no
\no
T^+ X^3&=&X^3 T^+ +q^{-2} \sqrt{1+q^2}X^+ \no
T^+ X^+&=&q^{-2} X^+T^+ \label{1.6} \\
T^+ X^-&=&q^2 X^- T^++q^{-1} \sqrt{1+q^2} X^3 \no
\no
T^- X^3&=&X^3 T^- +q \sqrt{1+q^2}X^- \no
T^- X^+&=&q^{-2} X^+ T^- +\sqrt{1+q^2}X^3 \label{1.7}\\
T^- X^-&=& q^2X^-T^-\nn
\ea
In the limit $q=1$ we obtain from relations (\ref{1.1})--(\ref{1.7}) 
the commutative $\b{R}^3$ with the Lie algebra $so(3)$ acting on it. 

As a consequence of the above relations it follows that there is
a central hermitean element, the $q$-deformed radius:
\ba
R^2&=&X^3 X^3-qX^+X^--\hq X^- X^+=q^2 \overline{X^3} X^3 +(1+q^{-2})
\overline{X^+} X^+, \\
\overline{R^2}&=&R^2. \nn
\ea
`Central' means that $R^2$ commutes with all the elements $X$ and $T$.

There is a well-known Casimir operator for the $su_q(2)$
algebra:
\be
\vec T^2=\frac{q^2}{\l^2} \tau^{\ha}+\frac{1}{\l^2}\tau^{-\ha}+
\tau^{-\ha}T^+T^--\frac{1+q^2}{\l^2}.                      \label{1.9}
\ee
We have introduced the group-like element
\be
\tau=1-\l T^3                              \label{1.10}
\ee
and the elements $\tau^{\ha}$ and $\tau^{-\ha}$ as an extension of the algebra.
We shall extend the algebra by the element $R=(R^2)^{\ha}$ and
$R^{-1}=(R^2)^{-\ha}$ as well.

The $\tau X$ and $\tau T$ commutation relations can be obtained from the $T^3X$ and
$T^3T$ relations and vice versa. They are
\ba\label{1.11}
\tau X^3&=&X^3 \tau \no
\tau X^+&=&q^{-4} X^+ \tau \\
\tau X^-&=&q^4 X^- \tau \nn
\ea
and
\ba
\tau T^3&=&T^3 \tau \no
\tau T^+&=&q^{-4} T^+ \tau \\
\tau T^-&=&q^4 T^- \tau. \nn
\ea

The definition of the orbital angular momentum as it was given in 
Ref. \cite{L.W.W.} can be best formulated in terms of the elements
\ba
L^+&=& {1 \over q^2 \sqrt{1+q^2}} \tau^{-\ha}T^+ \no
L^-&=& -{1 \over q^3 \sqrt{1+q^2}} \tau^{-\ha}T^- \label{1.14} \\
L^3&=& {1 \over q^2 (1-q^2)} \left( \tau^{-\ha} -1- {\l^2 \over 1+q^2}
\vec T^2\right) \nn
\ea
As the $q$-generalization of 
the fact that orbital angular momentum is orthogonal to the
coordinate vector we impose the constraint
\be
L \circ X=L^3 X^3-q L^+ X^--\hq L^- X^+=0.            \label{1.13}
\ee
We shall see that this defines orbital angular momentum uniquely.

\sect{The \boldmath{$t$} Algebra:}

The algebra introduced in the previous chapter allows a 
homomorphism of the $T$ algebra into the $X$ algebra. This was first
seen in Ref.\cite{C.F.M.}. We find this homomorphism by interpreting
Eqns. (\ref{1.5}), (\ref{1.6}) and (\ref{1.7}) as inhomogeneous
equations which can be solved for $T$ in terms of $X$. We first
construct a particular solution $t$ and exhibit the homomorphism
\ba
T^+ &\mapsto &t^+=-{\frac{1}{\l q^3}}\sqrt{1+q^2} X^+(X^3)^{-1}\no
T^- &\mapsto &t^-=\frac{q^2}{\l}\sqrt{1+q^2} X^-(X^3)^{-1} \label{2.1}\\
T^3 &\mapsto &t^3=\frac{1}{\l}\left(1+R^2(X^3)^{-2}\right). \nn
\ea
Here we extend the algebra by the inverse of $X^3$. 
To establish the homomorphism we have to use (\ref{1.1}) to show 
that the $t$ elements satisfy (\ref{1.3}). 
Furthermore the relations (\ref{1.5}) to (\ref{1.7}) are fullfilled 
by the $t$ elements. It is due to (\ref{1.2})
that they satisfy (\ref{1.4}) as well.

There are additional relations for the $t$ elements that follow from
(\ref{1.1}). They are
\be
\tau_t=1-\l t^3=-R^2(X^3)^{-2} \label{2.2}
\ee
and
\ba
t^+ t^-&=&-\frac{1}{\l^2}\left(1+q^2 \tau_t\right) \label{2.3}\\
t^- t^+&=&-\frac{1}{\l^2}\left(1+\nq \tau_t\right). \nn
\ea
It follows that the Casimir operator for  
the $t$ algebra takes a definite value and that in the notation of
Appendix \ref{appA}, where $\overline m_t$ and $d_t$ are defined,
\be
\vec T^2=-\frac{1+q^2}{\l^2}, \qquad \overline m_t=0, \qquad
d_t=-\frac{q^2}{\l}.                            \label{2.4}
\ee

This value of the Casimir operator and the sign of $\tau_t$, which is negative,
show that the $t$ algebra cannot be represented by the well-known
finite dimensional representations of the $T$ algebra
\cite{Wess}. In Appendix \ref{appA} we shall show that there are
infinite-dimensional representations of the $T$ algebra among which
there is one satisfying (\ref{2.2}), (\ref{2.3}) and (\ref{2.4}).
The representation is uniquely determined by these conditions,
we present it here:
\ba
t^3 \ket{m_t}&=&\frac{1}{\l}\left(1+q^2 q^{-4 m_t}\right) \ket{m_t} \no
t^+ \ket{m_t}&=&\frac{1}{\l q}\sqrt{q^{-4 m_t}-1} \ket{m_t+1} \label{2.5}\\
t^- \ket{m_t}&=&\frac{q}{\l}\sqrt{q^{-4 (m_t-1)}-1} \ket{m_t-1} \no
&&m_t \le 0. \nn
\ea
From (\ref{2.5}) it follows that
\be
t^+\ket{0}=0.
\ee
There is no state with positive $m_t$.

Eqns. (\ref{2.1}) allow us to express the elements $XR^{-1}$
in terms of the $t$ elements:
\ba
X^3 R^{-1}&=&\pm (-\tau_t)^{-\ha} \no
X^+ R^{-1}&=&\mp \frac{\l q^3}{\sqrt{1+q^2}}t^+(-\tau_t)^{-\ha} \\
X^- R^{-1}&=&\pm \frac{\l}{q^2\sqrt{1+q^2}}t^-(-\tau_t)^{-\ha} \nn
\ea
The two different signs are the signs of $X^3 R^{-1}=\sqrt{(X^3)^2 R^{-2}}$.
These elements can be viewed as homogeneous coordinates in the 
$\b{R}^3_q$ space.

The representations of these elements are now obtained from
(\ref{2.5}):
\ba
X^3 R^{-1}\ket{m_t}&=&\pm q^{2m_t-1} \ket{m_t} \no
X^+ R^{-1}\ket{m_t}&=&\mp \frac{q}{\sqrt{1+q^2}}\sqrt{1-q^{4m_t}}\ket{m_t+1}
\label{2.8}\\
X^- R^{-1}\ket{m_t}&=&\pm \frac{1}{\sqrt{1+q^2}}\sqrt{1-q^{4(m_t-1)}} 
\ket{m_t-1}
\nn
\ea
The different signs in  (\ref{2.8}) lead to
inequivalent irreducible representations of the $X$ algebra.

\sect{The \boldmath{$K$} Algebra}

We continue to consider the Eqns. (\ref{1.5}), (\ref{1.6}) and
(\ref{1.7}) as inhomogeneous equations that should be solved for
the $T$'s. We have found one particular solution (\ref{2.1})
and now move to the homogeneous part. This we do by the
Ansatz:
\ba
T^{\pm}&=&\Delta^{\pm}+ t^{\pm} \\
T^3&=& \Delta^3+t^3 \nn
\ea
Eqns (\ref{1.5}), (\ref{1.6}) and
(\ref{1.7}) become homogeneous equations for the $\Delta$'s.
\ba
X^\pm \Delta^+ &=&q^{\pm 2} \Delta^+X^\pm \\
X^3\Delta^+ &=&\Delta^+X^3 \no \no
X^\pm \Delta^- &=&q^{\pm 2} \Delta^- X^\pm \\
X^3\Delta^- &=&\Delta^- X^3 \no \no
X^\pm \Delta^3 &=&q^{\pm 4} \Delta^3 X^\pm \\
X^3\Delta^3 &=&\Delta^3 X^3 \nn
\ea
These equations suggest the further Ansatz:
\ba
K^\pm&=&\pm(-\tau_t)^{-\ha} \Delta^\pm,\\
K^3&=&(\tau_t)^{-1} \Delta^3. \nn
\ea
The element $\tau_t$ satisfies the relation 
(\ref{1.11}) and as a consequence all the $K$s 
commute with all the $X$s and therefore with all the $t$'s as well
\ba
K^A X^B &=&X^B K^A \\
K^A t^B &=&t^B K^A \nn
\ea

Now we turn to (\ref{1.3}) and compute the $KK$ relations:
\ba
q^{-1} K^+ K^--q K^-K^+&=&K^3\no
q^2 K^3 K^+-q^{-2} K^+ K^3&=&(q+q^{-1}) K^+ \\
-q^{-2}K^3 K^-+q^2 K^- K^3&=&(q+q^{-1}) K^- \nn
\ea
This is exactly the same algebra as (\ref{1.3}). Any
realization of the $K$-algebra will lead to a realization of the $T$-algebra:
\ba
T^{\pm}&=&t^{\pm}\pm(-\tau_t)^{\ha} K^{\pm} \label{3.8} \\
T^3&=&t^3 +\tau_t K^3 \nn
\ea
This is a relation which is familiar from the comultiplication of
two representations of the algebra (\ref{1.3}):
\ba
\Delta_{\beta}(T^3)&=&T^3 \otimes 1+\tau \otimes T^3 \\
\Delta_{\beta}(T^{\pm})&=&T^{\pm} \otimes 1\pm\sqrt{-\tau} \otimes T^{\pm}
\nn
\ea

This comultiplication will be discussed in Appendix \ref{appC}. It is 
adjusted to representations where the first factor has negative 
eigenvalues of $\tau$. We emphasize that the representations of the 
$t$ algebra in (\ref{3.8}) are restricted by the relations 
(\ref{2.3}) whereas for the $K$ algebra any representation
would do as long as we are not considering any conjugation
properties.

If we now demand the conjugation property (\ref{1.4}) for the $T$ algebra 
we find for the $K$ algebra:
\be
\overline{K^3}=K^3, \qquad \overline{K^+}=-\nq K^-.
\ee
Note the sign. The $K$ algebra belongs to the $SU_q(1,1)$
quantum group.

If we now use the condition (\ref{1.13}) for orbital angular momentum
 we will specify the $K$ algebra representation
uniquely as well. It needs some computation to express the
$L$ algebra (\ref{1.14}) in terms of the $t$ and $K$ algebras.
\ba
L^+&=& {1 \over q^2 \sqrt{1+q^2}} \left\{ (-\tau_t)^{-\ha}t^+ \otimes
(-\tau_k)^{-\ha}+1 \otimes (-\tau_k)^{-\ha} K^+ \right\}\no
L^-&=& -{1 \over q^3 \sqrt{1+q^2}} \left\{ (-\tau_t)^{-\ha}t^- \otimes
(-\tau_k)^{-\ha}-1 \otimes (-\tau_k)^{-\ha} K^- \right\} \\
L^3&=& {q^2-1 \over q^4 (q^2+1)} \left\{\frac{q^2}{\l^2}
(-\tau_t)^{\ha} \otimes (-\tau_k)^{\ha}+(-\tau_t)^{\ha} \otimes
(-\tau_k)^{-\ha} \left(-K^+ K^- + \frac{q^2}{\l^2}\right) \right.\no
&&\left. -\frac{1+q^2}{\l^2}
(-\tau_t)^{-\ha} \otimes (-\tau_k)^{-\ha} +t^- \otimes (-\tau_k)^{-\ha} K^+
-q^2 t^+ \otimes (-\tau_k)^{-\ha} K^-\right\}  \nn
\ea
This already shows that we should restrict the representations 
such that $\tau_k$ has negative eigenvalues. 
There is an additional reason for it. We shall see in Appendix~\ref{appC}
that the coproduct (\ref{3.8}) only leads to representations 
with positive eigenvalues of $\tau$ if $\tau_k$ has negative eigenvalues.
Only in this case the representations of $T$ can be
decomposed into finite-dimensional ones. We are here adding
this as an additional assumption - not knowing if it is really
necessary.

With this assumption it will follow from the comultiplication
rule of Appendix~\ref{appC} that we have to choose
\be
d_k=-\frac{1}{\l q^2}.                            \label{3.12}
\ee
Now we are ready to evaluate  (\ref{1.13}). This relation will be
true if and only if:
\be
\underline m_k=-1.                                         \label{3.13}
\ee
This is in the notation of Appendix \ref{appB}. For the Casimir operator
we find
\be
\vec T^2_k=-\frac{1+q^2}{\l^2}.                     \label{3.14}
\ee
This uniquely determines the $K$ algebra representation. The
generators of the orbital angular momentum will be denoted by
$T_{orb}$.

We find the result:
\ba
\T^3&=&t^3 \otimes 1+\tau_t \otimes K^3 \\
\T^{\pm}&=&t^{\pm} \otimes 1 \pm \sqrt{-\tau_t} \otimes K^{\pm}
\nn
\ea
where the $t$ and $K$ representations are determined by  (\ref{2.4})
and (\ref{3.12}), (\ref{3.13}).

We can add spin to orbital angular momentum:
\ba
T^3&=&\T^3 \otimes 1+\tau_{orb} \otimes S^3 \\
T^{\pm}&=&\T^{\pm} \otimes 1+\sqrt{\tau_{orb}} \otimes S^{\pm}
\nn
\ea
The spin operators $S$ can be in any finite-dimensional representations
of the $T$ algebra.

\sect{Representations of the \boldmath{$T_{orb}$} Algebra:}

The representation of the $K$ algebra that enters orbital angular
momentum is characterized by (\ref{3.12}), (\ref{3.13}) and
(\ref{3.14}):
\be
d_k= -\frac{1}{\l q^2}, \qquad \underline m_k=-1, \qquad \vec T^2_k=-
\frac{1+q^2}{\l^2}.
\ee
It is an infinite-dimensional representation with $m_k$ 
ranging from $0$ to $\infty$.
\ba
K^3 \ket{m_k}&=&\frac{1}{\l} \left(1+\nq q^{-4m_k}\right) \ket{m_k}\no
K^+ \ket{m_k}&=&\frac{1}{q\l} \sqrt{1-q^{-4(m_k+1)}} \ket{m_k+1} \label{4.2}\\
K^- \ket{m_k}&=&-\frac{q}{\l} \sqrt{1-q^{-4m_k}} \ket{m_k-1} \no
K^- \ket{0}&=&0, \qquad m_k\ge 0 \nn
\ea

The representation $T_{orb}$ of orbital angular momentum is the 
tensor product of this representation and the $t$ representation given
in (\ref{2.5}). The eigenstates of $\T^3$ are characterized by the
two numbers $m_t$ and $m_k$.
\ba
\T^3 \ket{m_t,m_k}&=&\frac{1}{\l} \left(1-q^{-4(m_t+m_k)}\right)\ket{m_t,m_k}
                                                    \label{4.3}\\
\T^+ \ket{m_t,m_k}&=&\frac{1}{\l q} \sqrt{q^{-4 m_t}-1} \ket{m_t+1,m_k}\no
&&+\frac{1}{\l} q^{-2m_t} \sqrt{1-q^{-4 (m_k+1)}} \ket{m_t,m_k+1}
\no
\T^- \ket{m_t,m_k}&=&\frac{q}{\l} \sqrt{q^{-4 (m_t-1)}-1} \ket{m_t-1,m_k}\no
&&+\frac{q^2}{\l} q^{-2m_t} \sqrt{1-q^{-4m_k}} \ket{m_t,m_k-1}
\nn
\ea
In this representation $\vec K^2$ and $\vec t^2$ are diagonal.
The transformation to the basis where $\vec T_{orb}^2$ is diagonal will
be constructed in the next chapter. The value of $\T^3$ in (\ref{4.3})
shows that we have found finite-dimensional representations of $T_{orb}$.

We obtain from (\ref{4.3})
\be
\tau_{orb} \ket{m_t,m_k}=q^{-4(m_t+m_k)}\ket{m_t,m_k}.   \label{4.4}
\ee

The representation of the $X$ algebra can be obtained from (\ref{2.8}).
The element $R$ is central, it will be diagonal in the $m_t,m_k$ basis.
We denote the eigenvalue of $R^2$ by
\be
R^2 \ket{m_t,m_k,M}=q^{4M+2}z_0^2 \ket{m_t,m_k,M}, \label{4.5}
\ee
where $z_0$ is an arbitrary parameter characterizing the radius.
Then we obtain from (\ref{2.8}) the representation of $X^3$:
\be
X^3 \ket{m_t,m_k,M}=q^{2(m_t+M)} z_0 \ket{m_t,m_k,M}. \label{4.6}
\ee
We have absorbed the sign in (\ref{2.8}) yielding inequivalent
representations in the sign of $z_0$ which is not determined by
(\ref{4.5}).

This and (\ref{4.4}) suggest that we should introduce a notation
characterizing the eigenvalue of $X^3$ by a quantum number as well
as the eigenvalue of $T^3_{orb}$.
\be
\nu=m_t+M, \qquad m=m_t+m_k                          \label{4.7}
\ee
In this notation we obtain the representation which was also found in 
\cite{Fiore1}:
\ba
X^3 \ket{M,\nu,m}&=&q^{2\nu} z_0 \ket{M,\nu,m} \no
R^2 \ket{M,\nu,m}&=&q^{4M+2}z_0^2 \ket{M,\nu,m} \no
T_{orb}^3 \ket{M,\nu,m}&=& \frac{1}{\l}\left(1-q^{-4m}\right) \ket{M,\nu,m}\no
X^+ \ket{M,\nu,m}&=&-\frac{q^2 z_0}{\sqrt{1+q^2}}\sqrt{q^{4M}-q^{4\nu}}
\ket{M,\nu+1,m+1} \label{4.8}\\
X^- \ket{M,\nu,m}&=&\frac{q z_0}{\sqrt{1+q^2}}\sqrt{q^{4M}-q^{4(\nu-1)}}
\ket{M,\nu-1,m-1} \no
T_{orb}^+ \ket{M,\nu,m}&=&\frac{1}{q^2-1} \sqrt{q^{4(M-\nu)}-1}
\ket{M,\nu+1,m+1}\no
&&+\frac{1}{\l}\sqrt{q^{4(M-\nu)}-q^{-4(m+1)}}\ket{M,\nu,m+1} \no
T_{orb}^- \ket{M,\nu,m}&=&\frac{q^2}{q^2-1} \sqrt{q^{4(M-\nu+1)}-1}
\ket{M,\nu-1,m-1}\no
&&+\frac{q^2}{\l}\sqrt{q^{4(M-\nu)}-q^{-4m}}\ket{M,\nu,m-1} \no
&&\nu \le M, m \ge \nu-M \nn
\ea

\sect{Reduction of the representation of \boldmath{$\T$}}

The above representation (\ref{4.3}) of $\T$ is a tensor product of two
representations, with $\vec t^2$ and $\vec K^2$ 
diagonal. We proceed with its decomposition into a sum of
irreducible representations characterized by the eigenvalues of $\vec \T^2$.
{From} the Appendix~\ref{appA} we know that for $d=\l d_t d_k=\l^{-1}$
the eigenvalues of $\vec T^2$ are $q[l][l+1]$.
Therefore we start with an Ansatz of the form
\ba
\ket{l,m}&=&\sum_{m_k,m_t} c_{l,m}^{m_k,m_t} \ket{m_t,m_k}, \qquad m_k \ge 0,
m_t \le 0 \\
\vec T_{orb}^2 \ket{l,m}&=&q [l][l+1] \ket{l,m}.
\label{ansatz}
\nn
\ea
According to (\ref{4.7}) $m=m_t+m_k$, so that we have
\be
c_{l,m}^{m_k,m_t}=c_{l,m}^{m_t} \delta_{m,m_t+m_k}.
\ee
{From} the definition (\ref{1.9}) of $\vec T^2$ and the Equations
(\ref{4.3}) we obtain a recursion relation for the 
coefficients $c_{l,m}^{m_t}$.
\ba
\label{t2}
\lefteqn{\left(q^{2l+2}+q^{-2l}-(q^2+1)q^{2(m+1)-4m_t}\right) c_{l,m}^{m_t}
=}  \nonumber\\
&&q^{2m+1} \left(\sqrt{(q^{-4m_t}-1)(q^{-4m_t}-q^{-4m})} c_{l,m}^{m_t+1}
\right.\\
&&\qquad\qquad\left.+\sqrt{(q^{4-4m_t}-1)(q^{4-4m_t}-q^{-4m})}
c_{l,m}^{m_t-1} \right)
\nonumber
\ea
A comparison with the $q$-difference Equation (\ref{difftp}) for the 
functions $\widetilde P^l_m$ defined in (\ref{defp}) and (\ref{defwp}) 
shows that (\ref{t2}) is solved by
\be
c_{l,m}^{m_t}=\left\{
\begin{array}{ll}
q^{m_t-m-1} \sqrt{1-q^{-2}} \widetilde P^{m}_l(\pm q^{2(m_t-1)-2m}) 
& \hbox{ for } m\ge 0 \\
q^{m_t-1} \sqrt{1-q^{-2}} \widetilde P^{|m|}_l(\pm q^{2(m_t-1)}) 
& \hbox{ for } m < 0
\end{array}
\right.
\ee
Note that $P^m_l$ is defined for $m \ge 0$ only.

The orthogonality condition (\ref{ortho}) for the functions
$\widetilde P^m_l$ suggests to start with the direct sum of two
representations of the form (\ref{4.3}), such that both signs of the
argument of $\widetilde P^m_l$ appear.
\ba
\ket{l,m}&=&\sum_{\sigma=\pm 1} \sum_{m_t} c_{l,m}^{m_t,\sigma} 
\ket{m_t,m_k,\sigma}
\label{transf} \\
c_{l,m}^{m_t,\sigma}&=&\left\{
\begin{array}{ll}
\sqrt{1-q^{-2}} q^{m_t-1-m}\widetilde P^{m}_l(\sigma q^{2(m_t-m-1}) 
& \hbox{ for } m\ge 0 \\ \\
\sqrt{1-q^{-2}} q^{m_t-1}\widetilde P^{|m|}_l(\sigma q^{2(m_t-1)}) 
& \hbox{ for } m < 0
\end{array}
\right.
\label{solution}
\ea

We know that $m_k \ge 0$, $m_t \le 0$ and $m=m_t+m_k$, thus
$m_t$ is restricted by $m_t \le 0$ and $m \ge m_t$.
The last condition comes into effect for negative values of $m$.
Note that if  $m_t$ takes its largest allowed value, the coefficient
of $c_{l,m}^{m_t+1}$ in (\ref{t2}) vanishes. We are free to choose
this $c$ to be zero. For $m \ge 0$ it then follows from (\ref{t2})
that $c_{l,m}^{m_t}=0$ for $m_t >0$ and for
$m<0$ the same is true for $m_t>m$.

The values of $l$ are restricted by the condition
$|m| \le l$, as seen from (\ref{condition}).
This is obviously consistent with the recursion formula (\ref{t2}).

We have chosen the normalization in (\ref{solution}) in such a way that 
according to (\ref{ortho}) the eigenfunctions of $\vec \T^2$
are orthonormal:
\be
\begin{array}{l}
q^{-1}\l \displaystyle{\sum_{\sigma=\pm 1} \sum_{m_t=-\infty}^{\min\{0,m\}}} 
q^{2(m_t-1)-m-|m|} \widetilde P^{|m|}_l(\sigma q^{2(m_t-1)-m-|m|})
\widetilde P^{|m|}_{l'}(\sigma q^{2(m_t-1)-m-|m|}) \\ \\
= \delta_{l,l'} 
\end{array}
\ee
To see this for $m<0$ it is enough to shift the summation variable 
$m_t \rightarrow m_t +m$.

We now assume that the two representations with $\sigma=+1$ and
$\sigma=-1$ also lead to a different sign of $z_0$ in (\ref{4.6}).
\be
X^3 \ket{m_t,m_k,M, \sigma}=q^{2 (m_t+M)} \sigma |z_0| \ket{m_t,m_k,M, \sigma}
\ee
Then it follows from (\ref{complete}) that the functions $\widetilde P^m_l$
satisfy the following completeness relation 
\be
\begin{array}{l}
q^{-1}\l  \displaystyle{\sum_{l=0}^{\infty}} q^{m_t+m'_t-2-m-|m|}
\widetilde P^{|m|}_l(\sigma q^{2(m_t-1)-m-|m|})
\widetilde P^{|m|}_{l}(\sigma' q^{2(m'_t-1)-m-|m|})\\ \\
=\delta_{\sigma,\sigma'}\delta_{m_t,m'_t}
\end{array}                            \label{compl}
\ee
This construction shows that for fixed $m$, $\vec \T^2$ is a 
selfadjoint operator in the basis $\ket{m_t,m_k,\sigma}$, $\sigma=\pm 1$, and
that the transformation from the basis $\ket{m_t,m_k,\sigma}$ to the basis 
$\ket{l,m}$ is an isometry.

\appendix
\section*{Appendices}

\sect{Representations of the \boldmath{$T$} Algebra}\label{appA}

When constructing representations of the $T$ algebra, we are aiming at 
representations where $T^3$ is selfadjoint (or essentially selfadjoint).
This allows us to assume $T^3$ to be diagonal:
\be
T^3 \ket{m}=f(m) \ket{m}.                   \label{A.1}
\ee
The eigenvalue of $T^3$ is $f(m)$, $m$ is a labelling of the
eigenstates.

The second equation of (\ref{1.3}) shows that $T^+\ket{m}$ is again 
an eigenstate of $T^3$, we choose the labelling such that this
state is labelled by $m+1$:
\be
T^+\ket{m}=c_m \ket{m+1}.                  \label{A.2}
\ee
The relation (\ref{1.3}) leads to a recursion formula for $f(m)$:
\be
f(m+1)=\frac{1}{q^4}f(m)+\nq (q+q^{-1}).
\ee
This recursion formula has the solution
\be
f(m)=\frac{1}{\l}-d q^{-4m}.                          \label{A.4}
\ee
{From} $\overline{T^3}=T^3$ follows that $dq^{-4m}$ has to be real.
We take $d$ and $m$ to be real.

For the operator $\tau$ of (\ref{1.10}) follows
\be
\tau \ket{m}=\lambda d q^{-4m} \ket{m}.                       \label{A.5}
\ee
From the conjugation properties of $T^+$ it follows that
\be
T^- \ket{m}=q^2 c^*_{m-1} \ket{m-1}.                  \label{A.6}
\ee
The third equation of (\ref{1.3}) is the conjugate of the second
one.
The first equation of (\ref{1.3}) amounts to a recursion formula
for $c^*_m c^{}_m$
\be
q c^*_{m-1} c^{}_{m-1}-q^3 c^*_m c^{}_m=f(m).               \label{A.7}
\ee
This recursion formula can be solved:
\be
c^*_m c^{}_m=\frac{1}{\l}\left\{-\frac{1}{q^2 \l}+\alpha \l q^{-2m}-
\frac{d}{q^4}q^{-4m}\right\}.                          \label{A.8}
\ee
The real parameter $\alpha$ is not determined by (\ref{A.7}).

We see that $c^*_m c^{}_m$ becomes negative for $m \rightarrow \infty$.
This is not allowed.
There has to be a largest value of $m$, say $\overline{m}$, such that
\be
c^*_{\overline{m}} c_{\overline{m}}=0.
\ee
Then it follows from (\ref{A.2}) that $T^+$ does not lead to a state with a
larger value than $\overline{m}$.
To analyze this situation we introduce the function:
\ba
x&=&q^{-2m} \label{A.10}                                            \\
h(x)&=&\left\{-\frac{1}{q(q^2-1)}+\alpha \l x-\frac{d}{q^4}x^2\right\}
\nn
\ea
The function $h(x)$ is negative for $x=0$, the sign of 
$h(x)$ for $x \rightarrow \infty$ depends on the sign of $d$.
In any case $h(x)$ has to have a zero for positive $x$ to represent $c^*c$.
We have to demand
\be
x_1=q^{-2 \overline m}, \qquad h(x_1)=0.            \label{A.11}
\ee
The parameter $\alpha$ can now be expressed in terms of $\overline m$.
\be
\label{A.12}
\alpha=\frac{1}{\l}\left\{\frac{1}{\l} q^{2(\overline m-1)}+d q^{-2(
\overline m+2)}\right\}.
\ee
If $\alpha$ takes this value $h(x)$ has the two zeros:
\be
x_1=q^{-2 \overline m}, \qquad x_2=\frac{1}{\l d}q^{2(\overline m+1)}.
                                                                \label{A.13}
\ee
We obtain for $c^* c$:
\be
c^*_m c^{}_m=-\frac{d}{\l q^4}\left(q^{-2m}-q^{-2\overline m}\right)
\left(q^{-2m}-\frac{1}{\l d} q^{2(\overline m+1)}\right).        \label{A.14}
\ee
The representation is characterized by two parameters, $\overline m$ and
$d$.
We use this parameter in the explicit form of the matrix elements:
\ba\label{A.15}
T^3 \ket{m}&=&(\frac{1}{\l}-d q^{-4m}) \ket{m} \no
T^+ \ket{m}&=&\nq \sqrt{\frac{d}{\l}(q^{-2m}-q^{-2\overline m})\left(
\frac{1}{\l d} q^{2(\overline m+1)}-q^{-2m}\right)} \ket{m+1} \\
T^- \ket{m}&=&\sqrt{\frac{d}{\l}(q^{-2(m-1)}-q^{-2\overline m})\left(
\frac{1}{\l d} q^{2(\overline m+1)}-q^{-2(m-1)}\right)} \ket{m-1} \no
\tau \ket{m}&=&d \l q^{-4m} \ket{m}\no
\vec{T}^2 &=&\frac{1}{\l^2 \sqrt{\l d}} \left(q^{2(\overline m+1)}+\l
d q^{-2 \overline m} \right)-\frac{1+q^2}{\l^2} \nn
\ea

Let us now have a closer look at the condition $c^*c \ge 0$. For this 
purpose we discuss the three cases $d>0$, $d=0$ and $d<0$
separately.

\medskip
\pn
$d>0$

There has to be a smallest value of m, say $\underline{m}$, such that
$\ket{\underline{m}}\neq 0 $ and $T^- \ket{\underline m}=0$,  therefore
\be
c^*_{\underline m-1} c^{}_{\underline m-1}=0.
\ee
{From} (\ref{A.14}) follows for $\overline m \ge 0$
\be
d=\frac{1}{\l}, \qquad \underline m=-\overline m.
\ee

The number of states between $\overline m$ and  $\underline m$ has to be
integer:
\be
2 \overline m+1 \equiv n.
\ee

This shows that $\overline m$ has to be integer or half integer and we found the
$2 l+1$-dimensional representation ($\overline m=l$) of $so_q(3)$.

\medskip
\pn
$d=0$

In this case $h(x)$ is a linear function:
\be
h(x)=-\frac{1}{q(q^2-1)}+\alpha \l x.
\ee
Now $\alpha$ has to be positive for $h$ to have a zero for positive $x$. 
{From} (\ref{A.12}) follows
\be
\alpha=\frac{1}{q^2 \l^2}q^{2\overline m}.
\ee
The representation can be obtained from (\ref{A.15}). The parameter
$\overline m$ that characterizes the representation can take any 
real value. The representation is infinite-dimensional, however,
$\tau$ is not invertible.

\medskip
\pn
$d<0$

This is the situation that arises for the $t$ algebra, as can be seen
{from} (\ref{2.2}). In this case $x_2$ is negative. We only have a 
largest value of $m$. The representation is 
infinite-dimensional and $\overline m$ is not restricted. The matrix elements
are obtained from (\ref{A.15}). We write them such as to exhibit the positive
square roots:
\ba
T^3 \ket{m}&=&(\frac{1}{\l}-d q^{-4m}) \ket{m} \no
T^+ \ket{m}&=&\nq \sqrt{-\frac{d}{\l}} \sqrt{(q^{-2m}-q^{-2\overline m})\left(q^{-2m}-
\frac{1}{\l d} q^{2(\overline m+1)}\right)} \ket{m+1} \\
T^- \ket{m}&=&\sqrt{-\frac{d}{\l}}\sqrt{(q^{-2(m-1)}-q^{-2\overline m})\left(
q^{-2(m-1)} -\frac{1}{\l d} q^{2(\overline m+1)}\right)} \ket{m-1} \no
\tau \ket{m}&=&d \l q^{-4m}\ket{m} \nn
\ea
$\tau$ has negative eigenvalues only, thus $\tau^{\ha}$ and $\vec T^2$ will
not be real. 

\sect{Representations of the \boldmath{$K$} Algebra}\label{appB}

The algebraic relations of the $K$ algebra are the same  as the relations
of the $T$ algebra, they are different only as a $*$ algebra:
\ba
\overline{K^3}=K^3, &\quad \overline{K^+}=-\displaystyle{\nq}
K^-  \label{B.1} \\
\overline{T^3}=T^3, &\quad \overline{T^+}=\displaystyle{\nq} T^- \nn
\ea
This makes the $K$ algebra a $su_q(1,1)$ algebra.

All the results that depend only on the algebraic relations are the
same as for the $T$ algebra.
\ba
\mbox{(\ref{A.1})}: & K^3 \ket{m}&\!\!=\;\phi(m) \ket{m} \\
\mbox{(\ref{A.2})}: & K^+ \ket{m}&\!\!=\;\gamma_m \ket{m} \\
\mbox{(\ref{A.4})}: & \phi(m)&\!\!=\;\frac{1}{\l}-d_k q^{-4m}
\ea
We again take $d_k$ and $m$ real.
\be
\mbox{(\ref{A.5})}:  \tau_k \ket{m}=\lambda d_k q^{-4m} \ket{m}
\ee
For $K^-$ there is a change in sign due to (\ref{B.1}):
\ba
\mbox{(\ref{A.6})}: & K^- \ket{m}&\!\!=\;-q^2 \gamma^*_{m-1} \ket{m-1} \\
\mbox{(\ref{A.8})}: & \gamma^*_m \gamma_m&\!\!=\;-\displaystyle{\frac{1}{\l}
\left\{-\frac{1}{q^2 \l}+\alpha \l q^{-2m}-\frac{d_k}{q^4} q^{-4m}\right\}} 
\ea
Now $\gamma^*_m \gamma_m$ becomes positive for $m \rightarrow \infty$,
we do not have to cut off the spectrum at a largest value of $m$.
We shall see that all the representations are infinite-dimensional.

We introduce the function $\kappa(x)$ analogous to $h(x)$ in (\ref{A.10}):
\be
\kappa(x)=\left\{ \frac{1}{q^2 \l^2} -\alpha x+\frac{d_k}{\l q^4}x^2\right\}, 
\qquad x=q^{-2m}. 
\ee
The representations of the $K$ algebra are:
\ba
K^3 \ket{m}&=&(\frac{1}{\l}-d_k q^{-4m}) \ket{m} \no
K^+ \ket{m}&=& \sqrt{\kappa(q^{-2m})} \ket{m+1} \\
K^- \ket{m}&=&-q^2\sqrt{\kappa(q^{-2(m-1)})} \ket{m-1} \nn
\ea
They are characterized by $\alpha$ and $d_k$ and restricted by the
condition $\kappa(q^{-2m}) \ge 0$.
To discuss this condition we determine the zeros of $\kappa(x)$
\ba
\kappa(x_{1,2})&=&0                                \label{B.10}\\
x_{1,2}&=&\frac{\l}{2 d_k q^{-4}} \left\{\alpha \pm \sqrt{\alpha^2-
4 d_k q^{-6} \l^{-3}} \right\}
\nn
\ea
We discuss the cases $d_k>0$, $d_k=0$  and $d_k<0$
separately and start with
\medskip
\pn
$d_k>0$:

In this case $\kappa(x)$ has no positive zero for 
$\alpha< \displaystyle{2q^{-3}\sqrt{d_k\lambda^{-3}}}=\alpha_0$.
The range of $m$ is not restricted, it can be of the form $m_0+n$,
$n \in \b{Z}_0$. If $\alpha \ge \alpha_0$ we will have two positive zeros
and $\kappa(x)$ can be written in the form
\be
\kappa(x)=\frac{1}{q^2 \l^2 x_1 x_2} (x-q^{-2 \overline{m}})
(x-q^{-2 \underline{m}}).
\ee
The values of the zeros $x_1$, $x_2$ determine the parameters 
$\alpha$ and $d_k$ and therefore the representation.
\ba
d_k&=&\frac{q^2}{\l}q^{2(\overline{m}+\underline{m})} \\
\alpha&=&\frac{1}{q^2 \l^2}(q^{2 \underline{m}}+q^{2 \overline{m}})
\nn
\ea
There are now two inequivalent representations. We find that in one
representation the allowed values of $m$ are
\be
m \le \overline{m}, \qquad m=\overline{m}, \overline{m}-1,
\overline{m}-2, \ldots  \label{B.13}
\ee
For the other representation we find
\be
\label{B.14}
m > \underline{m}, \qquad m=\underline{m}+1,\underline{m}+2, \ldots  
\ee
Now we consider
\medskip
\pn
$d_k=0$:

The function $\kappa$ becomes linear. It is positive at $x=0$ and, 
depending on $\alpha$, stays positive or becomes negative for
$x \rightarrow \infty$.
If $\alpha<0$ there is no restriction in the range of $m$, $m=m_0+n,$
$n \in \b{Z}$. If $\alpha>0$ there is a lowest eigenvalue of $m$, 
we are at the situation of (\ref{B.14}).

Finally we consider 
\medskip
\pn
$d_k<0$:

In this case $\kappa(x)$ is positive for $x \rightarrow 0$ and negative for
$x \rightarrow \infty$. There is one zero for $x>0$. This can also be
seen from (\ref{B.10})
because the square root will now be larger than $\alpha$. The
relevant zero of $\kappa(x)$ is:
\ba
x_1&=&\frac{\l}{2 |d_k| q^{-4}} \left\{-\alpha+\sqrt{ 
\alpha^2+4|d_k|q^{-6}\l^{-3}} \right\} \\
&=&q^{-2\underline m_k} \nn
\ea
Now all values of $\alpha$ are allowed. The range of $m$ will be as in
(\ref{B.14}).

For orbital angular momentum we encounter the representation
\be
d_k=-\frac{1}{q^2 \l}, \qquad x_1=q^2.
\ee
This leads to $\alpha=0$ and 
\be
\kappa(q^{-2m})=\frac{1}{\l^2q^2}\left(1-q^{-4(m+1)}\right).
\ee
The respective representation is shown in (\ref{4.2}).

\sect{Comultiplication}\label{appC}

The standard comultiplication rule for the algebra (\ref{1.3})
is:
\ba
\Delta(T^3)&=&T^3\otimes 1+\tau \otimes T^3 \label{C.1}\\
\Delta(T^{\pm})&=&T^{\pm} \otimes 1+\tau^{\ha} \otimes T^{\pm} \nn
\ea
As a consequence, $\tau$ is group-like:
\be
\Delta(\tau)=\tau \otimes \tau.                 \label{C.2}
\ee
The algebra (\ref{1.3}) is the same for the $T$ algebra
and the $K$ algebra, they are distinguished by their 
conjugation properties (\ref{B.1}).

As long as $\tau^{\ha}$ is hermitean, (\ref{C.1}) will respect the
conjugation properties  and we have a comultiplication within 
the $T$ algebra or the $K$ algebra respectively. From (\ref{C.2})
follows that $\Delta(\tau^{\ha})$ will be hermitean if $\tau^{\ha}$ is.
\be
\Delta(\tau^{\ha})=\tau^{\ha} \otimes \tau^{\ha}
\ee
If $\tau^{\ha}$ is not hermitean $\Delta(T)$ will  have no definite 
conjugation properties even if $T$ has.

We now turn to the product of representations as it follows from the
comultiplication rule (\ref{C.1}). If we have two representations 
of the algebra (\ref{1.3}) we obtain a new one by the rule
\ba
\Delta(T^3)&=&T_1^3\otimes 1+\tau_1 \otimes T_{2}^3 \label{C.4} \\
\Delta(T^{\pm})&=&T_1^{\pm} \otimes 1+\tau_1^{\ha} \otimes T_{2}^{\pm} \nn
\ea
{From} the discussion above follows that we can multiply two 
representations of the $T$ algebra ($K$ algebra) to obtain
a representation of the $T$ algebra ($K$ algebra)
as long as $\tau_1^{\frac{1}{2}}$ is hermitean. From now on we shall drop the indices $1$
and $2$ again, first and second representations will be defined by
the position in the product (\ref{C.4}).

That the $\tau^{\frac{1}{2}}$ of the first representation is hermitean means 
$d_1>0$. We shall discuss this situation first.
\medskip
\pn
$d_1>0$:

The product of two representations of the $T$ algebra ($K$ algebra)
will be a $T$ algebra ($K$ algebra). From (\ref{C.2}) follows
\be
d=\l d_1 d_2.                                        \label{C.5}
\ee
If $d_2$ is negative $d$ will be negative as well.

For the $T$ algebra $d$ positive restricts $d$ to be $d=\frac{1}{\l}$.
This characterizes the finite-dimensional representations.
{From} (\ref{C.5}) follows that the product of two
finite-dimensional representations is finite-dimensional as expected but also
that the product of a finite-dimensional representation
$(d_1=\frac{1}{\l})$ with an infinite-dimensional representation 
$(d_2<0)$ leads to $d<0$ and cannot be reduced to finite-dimensional 
representations.

For the $K$ algebra all representations are infinite-dimensional.

We now turn to the case that $d_1$ is negative, $\tau_1^{\frac{1}{2}}$ will 
be anti-hermitean.
\medskip
\pn
$d_1<0$:

In this case the product of two representations will in general
not have well-defined conjugation properties. We can, however, start 
from a modified comultiplication rule:
\ba
\Delta_{\beta}(T^3)&=&T^3\otimes 1+\tau \otimes T^3 \label{C.6}\\
\Delta_{\beta}(T^{\pm})&=&T^{\pm} \otimes 1 \pm (-\tau)^{\ha}
\otimes T^{\pm} \nn
\ea

If $(-\tau)^{\ha}$ is hermitean this rule allows us to multiply a
representation of the $T(K)$ algebra by a representation of the 
$K(T)$ algebra to obtain a $T(K)$ algebra.
\ba
&&T^3\otimes 1+\tau \otimes K^3 \\
&&T^{\pm} \otimes 1 \pm (-\tau)^{\ha} \otimes K^{\pm}
\nn
\ea
will be a representation of the $T$ algebra whereas
\ba
&&K^3\otimes 1+\tau_k \otimes T^3 \\
&&K^{\pm} \otimes 1 \pm (-\tau_k)^{\ha} \otimes T^{\pm}
\nn
\ea
will be a representation of the $K$ algebra.

For the comultiplication (\ref{C.6}) $\tau$ will be group-like as well
and it follows again that
\be
d=\l d_1 d_2.
\ee
But now $d_1$ is negative.

Of special interest is the case that $d_1$ and $d_2$ are both 
negative, then $d$ is positive. If we multiply $T \times K$ to obtain a 
$T$ algebra then we know that $d=\frac{1}{\l}$ and, as a consequence
\be
d_1 d_2=\frac{1}{\l^2}
\ee
to obtain a representation with well-defined conjugation properties.
This is exactly the case for the construction of the $T_{orb}$ algebra 
in the main part of this paper.

\sect{The big \boldmath{$q$}-Jacobi polynomials}\label{appD}

In this appendix we recall some basics about $q$-special functions
\cite{K.S.}, \cite{Koo}, \cite{Ga.Ra.}, in particular the big 
$q$-Jacobi polynomials.

First, we introduce some useful notation. 
The expressions
\be
[a] =\frac{q^a-q^{-a}}{q-q^{-1}} \stackrel{q \rightarrow 1} \rightarrow a,
\quad [a]!=\prod_{k=1}^a [k]  \stackrel{q \rightarrow 1} \rightarrow a! 
\ee
are known as symmetric $q$-numbers and symmetric $q$-factorials
respectively. The corresponding $q$-binomial coefficient is
\be
\left[ \begin{array}{l} n\\k \end{array} \right]=
\left\{ \begin{array}{ll} 
\displaystyle{\frac{[n]!}{[k]![n-k]!}} & \hbox{ for } n \ge k, \\ \\
0 & \hbox{ for $n<k$ or } n,k<0.
\end{array}
\right.                                        \label{binomial}
\ee
Of course $\left[ \begin{array}{l} n\\k \end{array} \right]\stackrel{q
\rightarrow 1} \rightarrow \left( \begin{array}{l} n\\k 
\end{array} \right)$. There are also ``unsymmetric'' counterparts of these 
objects: the basic $q$-number
\be
\frac{1-q^a}{1-q} \stackrel{q \rightarrow 1} \rightarrow a
\ee
and the $q$-shifted factorial (Pochammer-symbol)
\be
(a;q)_k=\prod_{n=0}^{k-1} (1-aq^n), \qquad 
(a_1,a_2,\ldots,a_i;q)_k=\prod_{m=1}^i (a_m;q)_k.
\ee
The Jackson integral of a function $f(x)$ is defined for $q>1$ by
\be
\int_0^a d_{q^{-1}} x \:f(x)=(1-q^{-1}) \sum_{\nu=0}^{\infty} aq^{-\nu} 
f(aq^{-\nu}) 
\label{jackson}
\ee
With the help of the $q$-shifted factorials, the basic hypergeometric
function can be introduced
\ba
\lefteqn{{}_r\phi_s\left(\left.
\begin{array}{l} a_1,\ldots,a_r \\ b_1,\ldots,b_s \end{array}
\right| q^{-1};x \right)=} \label{hyper}\\
&&\sum_{k=0}^{\infty}\frac{(a_1,\ldots,a_r;q^{-1})_k}
{(b_1,\ldots,b_s;q^{-1})_k}(-1)^{(1+s-r)k} q^{-\ha (1+s-r)k(k-1)} 
\frac{x^k}{(q^{-1};q^{-1})_k}
\nn
\ea
This series plays in the theory of $q$-special functions a role
analogous to that of the hypergeometric series in the theory of usual
special functions. We have considered a base $q^{-1}$ here, because
in this way ${}_r\phi_s$ is well-defined for $q>1$, which is the case we are
interested in here.

The big $q$-Jacobi polynomials \cite{K.S.}, \cite{Koo} are constructed in
terms of the basic hypergeometric series as
\be
P_l(x;a,b,c;q^{-1})={}_3\phi_2\left(\left.
\begin{array}{l} q^l,abq^{-(l+1)},x \\ aq^{-1},cq^{-1} \end{array}
\right| q^{-1};q^{-1} \right).                          \label{jacobi}
\ee
For the applications we consider in this paper, we are interested in the case
\ba
\lefteqn{P^m_l(x)\equiv P_{l-m}(x;q^{-2m},q^{-2m},-q^{-2m};q^{-2}),
\qquad m \ge 0}   \label{defp} \\
&&= \sum_{k=0}^{l-m} (-1)^k
\frac{q^{-k(m+1)} (x;q^{-2})_k}{(-q^{-2(m+1)};q^{-2})_k} 
\left[\begin{array}{l} l-m \\ k \end{array} \right]
\left[\begin{array}{l} l+m+k \\ k \end{array} \right]
\left[\begin{array}{l} m+k \\ k \end{array} \right]^{-1} \nn
\ea
Notice that the $P_l^m$ are polynomials of the order $l-m$ in $x$. Due to 
the factor $\left[\begin{array}{l} l-m \\ k \end{array} \right]$, 
which vanishes for $k>l-m$ according to the 
definition (\ref{binomial}). The sum in (\ref{hyper}) actually becomes finite. Due to the same factor the polynomials 
$P_l^m$ vanish, if the condition $m \le l$ is not satisfied:
\be
P^m_l(x)=0 \quad \hbox{for } l<m.                       \label{condition}
\ee
The further condition $m \ge 0$ is necessary for the polynomials $P^l_m$
to be well-defined, due to the factor $(q^{-2(m-1)};q^{-2})_k$ in the
denominator of the basic hypergeometric function which otherwise vanishes for
negative $m$.

Some of the first big $q$-Jacobi polynomials are
\be
\begin{array}{lll}
P^0_0(x)=1, & P^0_1(x)=x, \\ \\
P^0_2(x)=\displaystyle{\frac{1}{q[2]}}([3]x^2-q^{-2}),&
P^0_3(x)=\displaystyle{\frac{x}{q^5 [2]}}([5]q^2 x^2-[3]),\\ \\
P^1_1(x)=1, & P^1_2(x)=x, \\ \\
P^1_3(x)=\displaystyle{\frac{1}{q^5[4]}}(q^4[5]x^2-1). 
\end{array}
\ee

{From} Ref. \cite{K.S.} and \cite{Koo} we learn that the polynomials 
$P^m_l(x)$ satisfy a recurrence relation
\be
x q^m [2l+1] P^m_l(x)=q^l [l+m+1] P^m_{l+1}(x)
+q^{-l-1} [l-m] P^m_{l-1}(x)                         \label{recpol}
\ee
a $q$-difference equation
\ba
\lefteqn{\left( q^{-1-2m} (q^{2l+1}+q^{-2l-1}) x^2 -q^{-4(m+1)}(q^2+1)
\right) P^m_l(x)=} \label{diff}\\
&&q^{-2(2m+1)}(x^2-1)P^m_l(xq^{-2})+(x^2-q^{-4(m+1)}) P^m_l(xq^2)    \nn
\ea
and the orthonormality condition
\be
\int_{-q^{-2(m+1)}}^{q^{-2(m+1)}} d_{q^{-2}}x \: \: w_l^m(x)
w_{l'}^m(x) P^m_l(x)  P^m_{l'}(x) = \delta_{l,l'}.
\label{orthonorm}
\ee
Here, the weight function $w_l^m$
is defined by
\ba
w_l^m(x) &\equiv& \sqrt{\frac{\left(q^{-4(m+1)};q^{-4}
\right)_{\infty}}{\left(q^{-4},q^{-4(m+1)-2};q^{-4}\right)_{\infty}
\left(-q^{-2};q^{-2}\right)_{\infty}}} \sqrt{\frac{[2m+1]}{2[2l+1]}} 
\label{defw}\\
&&\times q^{-\frac{1}{2}(l^2+l+2lm-3 m^2+m+3)}
\sqrt{\left(x^2 q^{4m};q^{-4}\right)_m}
\sqrt{\frac{(q^{-2};q^{-2})_{l-m}}{(q^{-2(2m+1)};q^{-2})_{l-m}}} \nn
\ea
Actually, as it is done e.g. in \cite{Koo}, the big $q$-Jacobi polynomials
can be alternatively defined as those polynomials in $x$ which
are orthonormal with respect to the Jackson integral with the weight function 
$w_l^m(x)$ in (\ref{defw}).

The weight function has the following scaling properties
\ba
w_l^m(xq^{-2})&=&w_l^m(x)\sqrt{\frac{(1-x^2)}{(1-x^2q^{4m})}}\; , 
\label{scal} \\
w_{l-1}^m(x)&=&w_l^m(x)q^l\sqrt{\frac{[l+m][2l+1]}{[l-m][2l-1]}}\; . \nn
\ea
It is useful for the purposes of this paper to absorb the weight
function in the definition of the polynomials themselves and to introduce the
further notation
\be
\widetilde P_l^m(x)=w_l^m(x) P_l^m(x).                            \label{defwp}
\ee 
With the help of (\ref{scal}) it turns out that (\ref{recpol}) and 
(\ref{diff}) become respectively
\be
x \: q^{m+1} \widetilde P_l^m(x)= \sqrt{\frac{[l-m+1][l+m+1]}
{[2l+1][2l+3]}} \;\widetilde P_{l+1}^m(x) +\sqrt{\frac{[l+m][l-m]}
{[2l+1][2l-1]}} \;\widetilde P_{l-1}^m(x)                    \label{rectp}
\ee
and
\ba
\label{difftp}
\lefteqn{\left((q^{2l+1}+q^{-2l-1})q^{-1}x^2-(q^2+1)q^{-2(m+2)}\right) 
\widetilde P_l^m(x)=} \\
&&q^{-2(m+1)}\sqrt{(x^2-1)(x^2 q^{4m}-1)}\;\widetilde P_l^m(xq^{-2})\no
&&+\sqrt{(x^2-q^{-4(m+1)+1})(x^2-q^{-4})} \;\widetilde P_l^m(xq^2) \nonumber
\ea
By using (\ref{orthonorm}) and the definition of the Jackson integral 
(\ref{jackson}) we obtain the following 
orthonormality condition for the functions $\widetilde P^m_l(x)$
\be
(1-q^{-2}) \sum_{\sigma=\pm 1} \sum_{n=-\infty}^0 q^{2(n-m-1)} \widetilde P^m_l
(\sigma q^{2(n-m-1)}) \widetilde P^m_{l'}(\sigma q^{2(n-m-1)}) = \delta_{l,l'}.
\label{ortho}
\ee
Moreover, the functions $\widetilde P_l^m(x)$ have the property that they
transform under a parity transformation like
\be
\widetilde P^m_l(-x)=(-1)^{l-m} \widetilde P^m_l(x).
\ee

In the particular case $m=0$ the big $q$-Jacobi polynomials become the 
big $q$-Legendre polynomials, which in the limit $q \rightarrow 1$ yield
the usual Legendre polynomials. In the same limit from the polynomials 
$P_l^m(x)$ we recover the Jacobi polynomials with the
normalization $P^m_l(1)=1$.

\sect{Diagonalization of \boldmath{$X^3$}}\label{appE}

In this appendix we study the transformation which is inverse to the 
transformation (\ref{transf}), (\ref{solution}) constructed in Section 5. 
We show how the 
big $q$-Jacobi polynomials can be used to diagonalize $X^3$ in the basis 
where $\vec \T^2$, $\T^3$, $R^2$ are diagonal.

The representation where $\vec \T^2$, $\T^3$, $R^2$ are diagonal
can be found in \cite{C.W.}, \cite{doer}, \cite{S.Schr}
\ba
\vec T_{orb}^2\ket{M,l,m} &=& q [l][l+1]\ket{M,l,m} \no
X^3 \ket{M,l,m} &=& r_0 q^{2M+m} \left\{
\sqrt{\frac{[l+m+1][l-m+1]}{[2l+1][2l+3]}}
\ket{M,l+1,m} \right.\no
&&\qquad\qquad\left. +\sqrt{\frac{[l+m][l-m]}{[2l+1][2l-1]}} 
\ket{M,l-1,m} \right\} \\
X^+ \ket{M,l,m} &=& r_0 q^{2M+m} \left\{ q^{-l}
\sqrt{\frac{[l+m+1][l+m+2]}{[2][2l+1][2l+3]}}
\ket{M,l+1,m+1} \right.\no
&&\qquad\qquad\left.-q^{l+1}{\sqrt{\frac{[l-m][l-m-1]}{[2][2l+1][2l-1]}}} 
\ket{M,l-1,m+1} \right\} \label{5.1} \no
X^- \ket{M,l,m} &=& r_0 q^{2M+m} \left\{ q^l
\sqrt{\frac{[l-m+1][l-m+2]}{[2][2l+1][2l+3]}}
\ket{M,l+1,m-1} \right.\no
&&\qquad\qquad\left.-q^{-l-1}{\sqrt{\frac{[l+m][l+m-1]}{[2][2l+1][2l-1]}}} 
\ket{M,l-1,m-1} \right\} \nn
\ea
where
\be
0 \le l < \infty, \qquad -l \le m \le l.
\ee 
We make the following Ansatz for an eigenfunction of $X^3$
\be
X^3 \sum_{M,l,m} d_{M,l,m}\ket{M,l,m}=z \sum _{M,l,m} d_{M,l,m}
\ket{M,l,m}, 
\ee
with $z$ the corresponding eigenvalue.
By using (\ref{5.1}) we obtain a recursion relation for the coefficients 
$d_{M,l,m}$
\ba
z \: \: d_{M,l,m} &=&r_0 \frac{q^{2M+m}}{\sqrt{[2l+1]}}
\left\{\sqrt{\frac{[l-m+1][l+m+1]}{[2l+3]}} d_{M,l+1,m} \right. 
\label{c}\\
&& \qquad\qquad\qquad\left. +\sqrt{\frac{[l+m][l-m]}{[2l-1]}} d_{M,l-1,m} 
\right\}.
\nn
\ea
A comparison with the recursion relation (\ref{rectp})
for the functions $\widetilde P^m_l$ defined in (\ref{defwp}) in terms of the 
Jacobi polynomials shows that a solution of (\ref{c}) is
\be
d_{M,l,m}^{\nu,\sigma}=\left\{ \begin{array}{ll}
\sqrt{1-q^{-2}} q^{\nu-M-1-m} \widetilde P^m_l(\sigma q^{2(\nu-M-1-m)}) 
& \hbox{ for } m \ge 0\\
\sqrt{1-q^{-2}} q^{\nu-M-1} \widetilde P^{|m|}_l(\sigma q^{2(\nu-M-1)}) 
& \hbox{ for } m <0
\end{array}\right.
\label{pol}
\ee
where $z= \sigma r_0 q^{-1+2 \nu}$, $\sigma=\pm 1$.
By comparing with the form of the eigenvalues of $X^3$ (\ref{4.8}) 
we see that we have to restrict 
\be
\nu,M \in \b{Z},\quad \nu \le M, \quad m \ge \nu-M.
\label{eigen}
\ee
Notice that the argument of the functions would correspond to 
$x=\cos \theta=\frac{z}{r}$ classically, apart from
the $q$-factor $q^{-(m+|m|)} \rightarrow 1$ for $q\rightarrow 1$.

$X^3$ is a self-adjoint operator in this representation. This was shown 
in Ref. \cite{Weich}.

Now, the set of eigenfunctions of a self-adjoint operator is complete,
therefore we expect a completeness relation to hold for the eigenfunctions 
of $X^3$. In fact, (\ref{ortho}) can be interpreted in this way.
As the sum (\ref{ortho}) contains two sums, one
where the argument of $P^l_m$ is positive and one where it is negative,
we obtain a representation where the eigenvalues of $X^3$ can have
both signs, so that we automatically find the direct sum of two
representations of the type (\ref{4.8}). The normalization of the coefficients
$d_{M,l,m}$ in (\ref{pol}) has been chosen in such a way as to yield 
exactly (\ref{ortho}).

As the eigenfunctions of a selfadjoint operator corresponding to
different eigenvalues are orthogonal, since the normalization constant is
already fixed by (\ref{ortho}), we argue that the following relation 
holds
\be
(1-q^{-2}) \sum_{l=0}^{\infty} q^{\nu+\nu'-2}
\widetilde P^{|m|}_{l}(\sigma q^{2(\nu-1)})
\widetilde P^{|m|}_{l}(\sigma' q^{2(\nu'-1)}) 
=\delta_{\nu,\nu'} \delta_{\sigma,\sigma'} 
\label{complete} 
\ee
where $\sigma,\sigma'=\pm 1$ are the signs of the argument of the functions
and $\nu \le \min\{m,0\}$.
This is an interesting result for itself about the Jacobi polynomials.

\end{document}